\documentclass{article}

\usepackage{graphicx,amsmath,amssymb,amsthm}

 \numberwithin{equation}{section}

\begin{document}

{\Large \center \sc Judgment

}

\vspace{2mm}

{\center \sc Ruadhan O'Flanagan\footnote{Sloan-Swartz Centre for
Theoretical Neurobiology, Salk Institute, 10010 North Torrey Pines Road,
CA 92037. {\it oflanagan@salk.edu}}

}

\vspace{2mm}

\hrule

\vspace{1mm}

\hrule

\begin{abstract}
\noindent 
The concept of a judgment as a logical action which introduces new
information into a deductive system is examined. This leads to a way
of mathematically representing implication which is distinct from the
familiar material implication, according to which ``If $A$ then $B$'' is
considered to be equivalent to ``$B$ or not-$A$''. This leads, in turn,
to a resolution of the paradox of the raven.  
\end{abstract}

\hrule

\section{Introduction}

This note examines the relationship between the theory of quantification
of information and the concept of judgment, which is the action whereby
a proposition which provides information is given.

The information presented here is intended to be self-contained, but the
interested reader can consult references \cite{shannon,good,kullback}
for background information regarding information, \cite{jaynes,
poppernature} for probability, \cite{goodweigh,jaynesev,goodevidence}
for evidence, \cite{jasche,martin} for judgment and
\cite{hempel1,hempel2,goodraven,mackie,harre,rosenkrantz} for the
paradox of the raven, also called the raven paradox or the paradox of
confirmation.  On a first reading, the reader may find it advisable to
skip mathematical proofs which appear laborious; the text should remain
largely comprehensible.

Section \ref{judgment} introduces the concept of a judgment in
terms of its relation to probability. The central point
to understand is that a judgment involves a {\it change} in
what is thought to be true rather than merely a statement
that something is true, or an assertion. 

Considering judgments instead of statements makes it necessary to
take into account the distinction between the subject and the predicate.
For example, ``John is taller than Mary'' and ``Mary is shorter
than John'' express the same relation between John and Mary
when the sentences are interpreted as mere statements of fact.
It makes no difference whether John or Mary is considered to
be the subject, because one proposition is true if and only
if the other proposition is.

If they are interpreted as judgments, though, they produce different
changes: In the judgment that ``John is taller than Mary,'' John is the
subject and ``taller than Mary'' is the predicate. Mary is used only
for reference, so if we are incorporating this new information about
John into what we already know, then we will increase our estimate of
John's height without changing what we think about Mary. 

Conversely, in the proposition that ``Mary is shorter than John'',
Mary is understood to be the subject, so when this judgment
is passed, it is what is believed about Mary that changes.

Sections \ref{information} and \ref{evidence} introduce the quantities
of information and evidence in the context of judgment. The amount of
information given when a judgment occurs measures how much the judgment
can accomplish, that is, what changes in the probabilities of propositions
can be effected by the judgment.

Evidence accumulated in favour of a proposition inclines one towards
the judgment through which that proposition would be given; one can set
a threshold for evidence and pass the judgment that a proposition is
true if the evidence in favour of that proposition exceeds the threshold.
It is also shown in section \ref{evidence} that quantities of evidence
determine what changes should be made when a proposition of
probability zero is given.

Section \ref{implication} examines the various ways in which
it can be established through a judgment that one proposition implies
another. The statement that one proposition implies another is
unambiguous because it refers to a pre-existing state of affairs.
If one is to judge that one proposition implies another, however, one
must decide what changes to make, and there are several different ways
to change the relationship between the propositions so that,
after the changes have been made, the one proposition implies
the other.

In section \ref{implication}, the rules relating quantities of
information to counterfactuals are examined. In this context,
counterfactuals are propositions of the form ``If $A$ then $B$''
where $P(A)=0$. This leads to a consideration of what changes
can be effected when a proposition with probability one is given.

Finally, section \ref{raven} explains how distinguishing between judgments
and statements resolves the paradox of the raven. The paradox arises
when one regards ``All ravens are black'' and ``All non-black things
are non-ravens'' as equivalent, which would lead one to believe
that evidence favouring one proposition must also favour the
other. The observation of a black raven, then, should provide
evidence that ``All non-black things are non-ravens'', and
the observation of a blue sky, which is a non-black non-raven,
should provide evidence that ``All ravens are black'', which
is absurd.

The equivalence between the two propositions disappears when
one considers the judgments through which the propositions
are given rather than the statements which assert that the
propositions are true. 

The propositions have different subjects and different predicates:
the judgment that ``All ravens are black'' changes what we think about
ravens, while the judgment that ``All non-black things are non-ravens''
changes what we think about non-black things. Evidence inclines us towards
a judgment, not a statement; when we see that the evidence is sufficient
to persuade us, we change what we believe, rather than merely stating
what we already believed.

Section \ref{raven} makes this qualitative solution
quantitative by explicitly calculating the amounts of evidence
provided by the observation of a black raven in favour of
the respective propositions.

\section{Judgment} 
\label{judgment}

The word ``circumstance'' is used here to refer to a situation in which
certain propositions are given and every proposition can be assigned a
probability. A judgment is a logical action which changes the circumstance
by giving a new proposition. If $A$ is judged to be true then, after the
judgment, the new value of $P(B)$ is equal to the old value of $P(B|A)$.

Since the probabilities of propositions depend on the circumstance,
an unambiguous notation would indicate, for each probability mentioned,
what the relevant circumstance was. $P_C(A)$ could be used, for example,
to denote the probability of $A$ in the circumstance $C$. Since we will
mostly be considering individual judgments here, though, we will drop
the $C$ and refer to the values of probabilities before the judgment as
the old values and the values after the judgment as the new values.

The expression $P(B|A)$ is read as ``the probability of $B$ given $A$''
but in this case $A$ is only hypothetically given; it is not actually
taken as an axiom in the present circumstance. When the judgment occurs,
$A$ is actually given.  $P(B|A)$ is more exactly described by the phrase
``the probability that $B$ would acquire if $A$ were to be given''.

When $P(A)\neq0$, $P(B|A)=P(BA)/P(A)$, where $BA$ has been used to
denote the proposition $B\mbox{\it\ and }A$. This is not a definition of
$P(B|A)$. Since it has already been said that $P(B|A)$ is the probability
that $B$ would acquire if $A$ were to be given, we are not free to define
$P(B|A)$ again. That $P(B|A)=P(BA)/P(A)$ when $P(A)\neq 0$ is something
which must be justified, or proven, instead of chosen. The justification
is easily understood when the equation is expressed in terms
of information, as it is in the next section.

The interpretation of probabilities is not important here. It is only
necessary to observe that $P(A)=1$ may provisionally be considered
as a sufficient proof of the proposition $A$, but that it is not a
logical impossibility for $P(A)$ to equal one and for $A$ to be false.
For example, if a real number, $x$, is chosen randomly from between zero
and one, so that the probability that $x$ lies in the interval $[c,d)$ is
$d-c$, then before the value of $x$ is revealed to us, the probability that
$x=1/2$ is zero, but it is not logically impossible that $x=1/2$. There
is in any case some number which $x$ equals, and the probability that $x$
equals that number is zero.

The inference from $P(A)=1$ to $A$ is therefore weaker than 
a perfect proof, and if $P(A)=1$ in a particular circumstance,
then it remains a logical possibility that $A$ is false. It
should therefore be possible to incorporate the information
that $A$ is false without encountering a logical contradiction,
although after this information has been incorporated, the
probability of $A$ will be zero instead of one. 

It will be said that a proposition, $A$, is ``believed'' when, in a
particular circumstance, $P(A)=1$, and, correspondingly, that it is not
believed when $P(A)\neq 1$. This word is used because it is weaker
than, for example, the word ``known'', which would suggest not
only that the proposition is believed, but also that it is true.

When $P(A)=0$, $P(B|A)$ cannot be determined by using the expression
$P(BA)/P(A)$. However, all that is necessary for the judgment to occur
is that every proposition $B$ be assigned a probability consistent with
the rules of probability\footnote{The rules relating probability to
propositions are: $P(A\mbox{\it\ or\ }B)=P(A)+P(B)-P(AB)$, $P(AB)+P(A\bar
B)=P(A)$ and $0\leq P(A) \leq 1$, where $\bar B$ denotes the negation
of $B$.}, which is enough to specify the new circumstance.

For example, the assignment:
$$
P(x=y|x=1/2)
=
\left\{
\begin{array}{lr}
0 & \mbox{\ \ \ if\ \ } y\neq \frac{1}{2} \\
1 & \mbox{\ \ \ if\ \ } y= \frac{1}{2} \\
\end{array}
\right.
$$
is consistent even though the expression
$P(x=y\mbox{\it\ and\ }x=\frac{1}{2})/P(x=\frac{1}{2})$ yields $0/0$.

If $P(A)=0$, then, the value of $P(B|A)$ is not constrained by the values
of probabilities of propositions in the current circumstance, while
if $P(A)>0$, it is. The inference from $P(A)=1$ to $A$ is therefore
legitimate insofar as $P(A)$ is guaranteed to equal one in every future
circumstance unless a proposition of probability zero is given, after
which $P(A)$ may take any value.

\section{Information}
\label{information}

\subsection{Introduction}

$i(B)\equiv -\log P(B)$ is the amount of information which would be
provided by the revelation that $B$ is true, since it is zero if $B$
is already believed, one bit if $P(B)=1/2$, two bits if $P(B)=1/4$ and
so on. 

If we are willing to say that information is something which can be
believed, in addition to saying that propositions can be believed, then
$i(B)$ can be called the amount of information which it is necessary to
believe in order to believe $B$. If $B$ implies and is implied by $CD$
and $C$ and $D$ are independent, that is, $P(CD)=P(C)P(D)$, then in order
to believe $B$ it is necessary to believe $C$ and also to believe $D$,
and the amount of information which it is necessary to believe in order
to believe $B$ is equal to the sum of the amount necessary to believe $C$
and the amount necessary to believe $D$.

It can be seen that $i(B)$ is the minimum amount of information
which must be given before $P(B)$ can reach one from the fact that:
$$P(B|A)=1 \Rightarrow 1=P(BA)/P(A)\leq P(B)/P(A) \Rightarrow i(A)\geq i(B)$$

If $A$ implies $B$ but is no less probable than it\footnote{$A$
cannot be more probable than $B$ if $A$ implies $B$, so if $A$ is no
more informative than $B$ and implies it then $P(A)=P(B)$.}, then $B$
also implies $A$, because: $$1=P(B|A)=P(BA)/P(A)=P(BA)/P(B)=P(A|B)$$ 

This applies only if $P(A)\neq0$. We can therefore infer from the truth
of the consequence of a proposition to the truth of the proposition
itself if the proposition is no less probable, or equivalently, no more
informative than the consequence.

\subsection{Common Information}

The amount of information which it is necessary to believe in
order to believe a proposition can never be negative, because
$P(A)\leq 1 \Rightarrow i(A)\geq 0$.

If we consider the quantity $i(A;B)\equiv i(A)+i(B)-i(AB)$, then
this quantity can be positive or negative.

When $i(A;B)$ is positive, less information is required to believe
both $A$ and $B$ than the sum of the amounts of information required to
believe $A$ and $B$. This is the condition under which three independent
propositions, $C$, $D$ and $E$ can be known of such that $A$ implies
and is implied by $CD$ and $B$ implies and is implied by $DE$, with
$i(D)=i(A;B)$. In that case, $i(A;B)$ could be called the amount
of information common to $A$ and $B$.

On the other hand, if the quantity $i(A;B)$ is negative, then in order
to believe both $A$ and $B$ it is necessary to believe more than the
sum of the amounts necessary to believe $A$ and $B$ individually. When
$i(A;B)$ is negative, then, $AB$ implies something additional besides
what is implied by $A$ and what is implied by $B$.  The question which
naturally arises is whether a proposition, $C$, can be known of, such
that $C$ is independent of $A$ and independent of $B$ but is implied by
$AB$, with $P(AB|C)=P(A)P(B)$, which in this case would be equivalent
to $i(C)=-i(A;B)$.

For a proposition such as $C$ to stand in that relationship to $A$ and
$B$, it would be necessary for $P(A)+P(B)$ to be less than or equal to
one. This is because if $C$ is independent of $A$ and $B$, then it must be
possible to learn that $C$ is false without changing the probabilities
of $A$ and $B$, so $P(A)+P(B)=P\left(A|\bar C\right)+P\left(B|\bar
C\right)$. If $C$ is discovered to be false, however, then $AB$ will be
believed to be false and hence the truth of $A$ would imply the falsity
of $B$ and vice versa, $P\left(AB|\bar C\right)=0$. That is, $A$
and $B$ would be mutually exclusive, and therefore $P\left(A|\bar
C\right)+P\left(B|\bar C\right)$ could not exceed one, and neither
could $P(A)+P(B)$.

Provided $i(A;B)<0$ and $P(A)+P(B)\leq 1$, though, a proposition such as
$C$, which $A$ and $B$ jointly imply but are individually independent of,
can be known of. $C$ could be called the independent consequence of $A$
and $B$.

It can be seen that:
$$
\begin{array}{rl}
i(A;B)<0 &\Rightarrow  \log\frac{P(AB)}{P(A)P(B)}<0 \\
& \Rightarrow P(AB)<P(A)P(B) \\
& \Rightarrow 1-P(A)-P(B)+P(AB)<(1-P(A))(1-P(B)) \\
& \Rightarrow P\left(\bar A\bar B\right) < P\left(\bar A\right)P\left(\bar B\right) \\
& \Rightarrow \log\frac{P\left(\bar A\bar B\right)}{P\left(\bar A\right)P\left(\bar B\right)}<0 \\
& \Rightarrow i\left(\bar A;\bar B\right)<0
\end{array}
$$
where it has been assumed in the second-last step that $P\left(\bar
A\right)>0$ and $P\left(\bar B\right)>0$. Now either $P(A)+P(B)\leq 1$
or $P\left(\bar A\right)+P\left(\bar B\right)\leq 1$ or both. Combining
this with the result above proves that either $A$ and $B$ or $\bar A$
and $\bar B$ will satisfy the conditions necessary to have an independent
consequence.

More generally, for any two propositions, $A$ and $B$, which are not
independent and with probabilities strictly between zero and one, two
of the quantities $i(A;B)$, $i\left(\bar A;B\right)$, $i\left(A;\bar
B\right)$ and $i\left(\bar A;\bar B\right)$ will be negative and two
will be positive. Of the two pairs of propositions which yield
negative quantities, it will be possible for at least one pair to have
an independent consequence.

\subsection{The Justification of $P(A|B)=P(AB)/P(B)$}

We can denote the amount of information which it is necessary to believe
in order to believe $A$ assuming that $B$ is true as $i(A|B)\equiv-\log
P(A|B)$. Whenever $P(B)\neq 0$, $i(A|B)$ is equal to $i(AB)-i(B)$.
This is equivalent to the statement that $P(A|B)=P(AB)/P(B)$ but is
more obviously true. The following three expressions clearly all refer
to the same quantity:
\begin{itemize}
\item
The amount of information which it is necessary to believe in order to
believe that $A$ is true if we can assume that $B$ is true
\item
The amount of information which it is necessary to believe in order to
believe that both $A$ and $B$ are true if we can assume that $B$ is true
\item
The amount of information which it is necessary to believe in order
to believe that both $A$ and $B$ are true minus the amount necessary
to believe $B$ (which has been assumed and therefore does not need to
be believed)
\end{itemize}

This can be taken as a justification of the statement that
$P(A|B)=P(AB)/P(B)$ whenever $P(B)\neq 0$. $i(A|B)$ is also equal to
$i(A)-i(A;B)$.

\subsection{How Much a Judgment can Accomplish}

For any two propositions, $A$ and $B$, such that $P(B)\neq 0$, $i(A|B)$ is
always greater than or equal to $i(A)-i(B)$ because $\log P(AB)/P(B)\leq
\log P(A)/P(B)$. If, when $B$ is given, $i(A)$ decreases to $i(A)-i(B)$,
then the judgment through which $B$ is given can not alter the probability
of any proposition which is independent of $A$ both before and after the
judgment. This is because, if $D$ is such a proposition and $P(D|B)\neq
P(D)$, then either $i(D|B)<i(D)$ or $i(\bar D|B)<i(\bar D)$, so either
$$i(AD|B)=i(A|B)+i(D|B)=i(A)-i(B)+i(D|B)<i(AD)-i(B)$$ 
or
$$i(A\bar D|B)=i(A|B)+i(\bar D|B)=i(A)-i(B)+i(\bar D|B)<i(A\bar D)-i(B)$$
either of which would contradict the statement that the amount of
information necessary to believe a proposition can decrease by at most
$i(B)$ when $B$ is given. In this sense, if all of the information given
in a judgment is contained within the information required to believe a
particular proposition, then that judgment can accomplish nothing else,
that is, no changes in the probabilities of propositions which are
independent of that particular proposition.

\section{Evidence}

\label{evidence}

\subsection{Quantifying Evidence}

If a coin is biased so that it lands on one particular side 90\% of the
time, but it is not known which side the coin is biased in favour of,
then repetitively tossing the coin and observing which side it most
frequently lands on provides a way to discover which side the bias
favours.

If $H$ is the hypothesis that the bias in favour of heads, and
$B_n$ is the proposition that the result of the $n^{\mbox{th}}$
toss is heads, then:
$$
\begin{array}{rl}

\log \frac{P(H|B_1B_2)}{P\left(\bar H|B_1B_2\right)}
& =\log \frac{P(B_1B_2|H)}{P(B_1B_2)}\frac{P(B_1B_2)}{P\left(B_1B_2|\bar H\right)}
\frac{P(H)}{P\left(\bar H\right)} \\
& =\log \frac{P(H)}{P\left(\bar H\right)}+\log \frac{P(B_1B_2|H)}{P\left(B_1B_2|\bar H\right)}
\end{array}
$$

If we can assume that $H$ is true, then observing the result of one coin
toss does not change the probability that the next toss will result
in heads; it remains at 90\%. Similarly, if we assume that $H$ is
false, then the probability that the second result will be heads
is 10\%, regardless of what the result of the first toss is. The
probabilities $P(B_1B_2|H)$ and $P\left(B_1B_2|\bar H\right)$
therefore factorize:
$$
P(B_1B_2|H)=P(B_1|H)P(B_2|H) \mbox{\ \ \ \ and \ \ \ \ }
P(B_1B_2|\bar H)=P(B_1|\bar H)P(B_2|\bar H)
$$
leading to:
$$
\log \frac{P(H|B_1B_2)}{P\left(\bar H|B_1B_2\right)}
=\log \frac{P(H)}{P\left(\bar H\right)} + \log \frac{P(B_1|H)}{P\left(B_1|\bar H\right)} + \log \frac{P(B_2|H)}{P\left(B_2|\bar H\right)}
$$

In this way, each occurrence of a heads provides an additional
contribution to the logarithm of the odds of $H$, independent
of the contributions provided by the results of the other coin
tosses.

The quantity which accumulates as more occurrences of heads are observed,
with each occurrence providing a contribution independent of the others,
and such that the accumulation of this quantity inclines one towards
the belief that the bias is in favour of heads, is colloquially called
evidence. We will use the notation $e(B_1\rightarrow H)$ to refer
to $\log \frac{P(B_1|H)}{P\left(B_1|\bar H\right)}$, and call it
the amount of evidence provided by the proposition $B_1$ in
favour of the proposition $H$.

No finite number of coin tosses can ever provide enough evidence to prove
that the bias is in favour of heads or of tails, because the logarithm of
the odds of a proposition must reach infinity in order for the probability
of that proposition to reach one. Instead one can introduce a threshold
for evidence, and if the amount of accumulated evidence in favour of a
proposition exceeds that threshold, one can pass the judgment that
the proposition is true, having been persuaded by the evidence.

The logarithm of the odds of $H$, $\log \frac{P(H)}{P\left(\bar H\right)}$,
could be called the amount of evidence in favour of $H$, and denoted by
$e(H)$. In cases where $P(H)=0$, or $P(H)=1$, however, this quantity
is infinite and does not change when new evidence in favour of $H$ or
against $H$ is given.  Keeping track of the amount of evidence which
has accumulated is therefore not accomplished by keeping track of $P(H)$.

If we want, for example, to be able to initially believe that the bias
of the coin is in favour of tails, $P(H)=0$, but change and believe
that it is in favour of heads when the accumulated evidence has reached
some threshold\footnote{Since a judgment which is made on the basis of
a finite amount of evidence always carries with it the risk of error,
it is wise to retain the possibility of making the opposite judgment
if enough evidence subsequently accumulates favouring the opposite
conclusion.}, then it would be necessary to keep track of the
quantity of accumulated evidence, $e_a(H)$, instead of $e(H)$. $e_a(H)$
would start at zero and receive additive contributions from each
observation of a coin toss. Until the threshold for $e_a(H)$ is reached,
$P(H)$ would remain at zero and $e(H)$ would remain at $-\infty$.

Evidence is related to information as follows:
$$
e(B\rightarrow A)=\log \frac{P(B|A)}{P\left(B|\bar A\right)}
=\log \frac{P(B|A)}{P(B)}+\log \frac{P(B)}{P\left(B|\bar A\right)}
=i(A;B)-i(\bar A;B)
$$

\subsection{When a Proposition with Probability Zero is Given}

The quantity
$$
e(B\rightarrow A)=\log \frac{P(B|A)}{P(B)}
+
\log \frac{P(B)}{P\left(B|\bar A\right)}
=i\left(B|\bar A\right)-i(B|A)
$$
has another significance with respect to judgment in circumstances
in which $P(A)=0$. In those circumstances, $P(B)=P\left(B|\bar A\right)$
so $i(B)=i\left(B|\bar A\right)$. 

In order to pass the judgment that $A$ is true, it is necessary
to know what changes to make to the probabilities of other
propositions such as $B$. This is equivalent to specifying
the changes that should be made to quantities of information
such as $i(B)$. If $i(B)=i\left(B|\bar A\right)$ initially,
then the quantity which must be subtracted from it when $A$ is given
is $i\left(B|\bar A\right)-i(B|A)=e(B\rightarrow A)$.

The amount of evidence provided by $B$ in favour of $A$ is therefore
equal to the amount by which $i(B)$ should change if $A$ is given in
a circumstance in which $P(A)=0$. 

In the example where the results of coin tosses are observed, it might
be initially believed that the bias is in favour of tails, $P(H)=0$,
in which case the probability that the next toss will land on heads,
$P(B_n)$, would be $P\left(B_n|\bar H\right)=0.1$.

If $H$ is given, $i(B_n)$ will change to: 
$$ 
i(B_n|H)=-\log 0.1 - e(B_n\rightarrow H)=-\log 0.1- \log \frac{0.9}{0.1} = -\log 0.9
$$
so $P(B_n)$ would change to 0.9, which is $P(B_n|H)$.

It is therefore sufficient to know the values of quantities of evidence
such as $e(B_n\rightarrow H)$ in order to know what changes to make to
the probabilities of propositions when a proposition of probability zero,
such as $H$, is given.

\subsection{Evidence and Independence}

If $P(A)=0$, then $i(B|A)=i(B)-e(B\rightarrow A)$ so if $A$ is
judged to be true, $e(B\rightarrow A)$ will be subtracted from
$i(B)$.

To undo the judgment through which $A$ was given, and return to the
original circumstance, the same quantity which was subtracted from $i(B)$
can be added to the new value of $i(B)$ to obtain the original value
of $i(B)$. The quantity which is to be added is then equal to the old
value of $e(B\rightarrow A)$, that is, the value which $e(B\rightarrow
A)$ had before $A$ was given.

Performing this addition must be equivalent to subtracting the new
value of $e\left(B\rightarrow \bar A\right)$ since this is what must
be subtracted from $i(B)$ when $P\left(\bar A\right)=0$ and $\bar A$
is given. Since $e(B\rightarrow A)$ is equal to $-e\left(B\rightarrow
\bar A\right)$ in every circumstance, this shows that $e(B\rightarrow
A)$ does not change when either $A$ or $\bar A$ is given in a circumstance
where its initial probability is zero, unlike, for example, $i(A;B)$.

Whenever $i(A;B)>0$, $i(\bar A;B)\leq 0$, and, correspondingly,
$i(A;B)<0$ implies that $i(\bar A;B)\geq 0$, so if the amount
of evidence that $B$ provides in favour of $A$ is zero, then
the amount of information common to $A$ and $B$ is zero.

The converse is not true, however. That is, it is possible for
$i(A;B)$ to be zero while $e(B\rightarrow A)$ is non-zero, namely
if $A$ has probability one. The condition that
$P(A)P(B)=P(AB)$ is therefore a weaker type of independence
than the condition that $e(B\rightarrow A)=0$. 

Even if $e(B\rightarrow A)=0$, it is still possible for
$e(A\rightarrow B)$ to be different from zero, namely if
$P(B)=1$. The quantity:
$$
e_m(A;B)\equiv i(A;B)-i(\bar A;B)-i(A;\bar B)+i(\bar A;\bar B)
=e(B\rightarrow A)-e(\bar B\rightarrow A)
$$
is what must equal zero in order to ensure that $A$ and
$B$ are fully independent, that is, that each of the
four quantities $i(A;B)$, $i(\bar A;B)$, $i(A;\bar B)$ and $i(\bar
A;\bar B)$ vanishes. $e_m(A;B)$ is the amount by which the
log of the odds of $A$ changes if $B$ is given in a circumstance
where $P(B)=0$. Since it is symmetric, it is also the amount
by which the log of the odds of $B$ changes if $A$ changes
from being believed to be false to being believed to be true.

$e_m(A;B)$ can be called the amount of evidence which $A$ and $B$
mutually provide in favour of one another, or the mutual evidence, since
it is a quantity of evidence which is symmetric under interchange of the
propositions $A$ and $B$, and it pertains to each proposition in respect
of the other. $i(A;B)$ is called common instead of mutual because the
information referred to pertains to each proposition as well as the other,
rather than in respect of the other. For reasons analogous to the case of
$e(B\rightarrow A)$, $e_m(A;B)$ does not change when any of $A$, $\bar
A$, $B$ or $\bar B$ is believed to be false but then given.

\section{Implication}

\label{implication}

In propositional logic it is customary to use the proposition
$B$ {\it or }$\bar A$ as the logical or mathematical representation
of the implication which is expressed in English as ``If $A$ then $B$''.

We can consider the information common to $\bar A$ and $B$ {\it or }$\bar A$:
$$
\begin{array}{rl}
i\left(\bar A;B\mbox{\it\ or }\bar A\right) & =
\log \frac{P\left(\bar A\mbox{\it\ and }(B\mbox{\it\ or }\bar A)\right)}{P\left(\bar A\right)P\left(B\mbox{\it\ or }\bar A\right)}
=
\log \frac{P\left(\bar A\right)}{P\left(\bar A\right)P\left(B\mbox{\it\ or }\bar A\right)} \\
& =
\log \frac{1}{P\left(B\mbox{\it\ or }\bar A\right)}
=i\left(B\mbox{\it\ or }\bar A\right)
\end{array}
$$

Unless it is already believed that $B$ is true or $A$ is false, this
quantity is positive.  This implies that:
$$
i\left(\bar A|B\mbox{\it\ or }\bar A\right)=i(\bar A)-i\left(B\mbox{\it\ or }\bar A\right)
< i(\bar A)
$$
which in turn implies that $P\left(\bar A|B\mbox{\it\ or }\bar A\right)>P\left(\bar A\right)$, or equivalently $P\left(A|B\mbox{\it\ or }\bar A\right)<P(A)$.

This result can be stated in words by saying that if it is judged that
$B$ is true or $A$ is false then the probability of $A$ must
decrease unless it was already believed that $B$ is true or $A$ is false.
Similarly, the probability of $B$ must increase.

If we examine the senses in which the English expression ``If $A$ then $B$''
can be used, however, we find that it is not always the case that
incorporating the information provided leads us to regard $A$ as
less probable and $B$ as more probable. Four cases can be distinguished:

\begin{itemize}
\item
Both $P(A)$ and $P(B)$ change.
\item
$P(A)$ remains unchanged; $P(B)$ changes.
\item
$P(A)$ changes; $P(B)$ remains unchanged.
\item
Both $P(A)$ and $P(B)$ remain unchanged.
\end{itemize}

We can examine each of the four cases in turn, using the letter $T$
to denote the proposition which is given by the judgment establishing
that $A$ implies $B$.

\vspace{5mm}

\subsection{The Judgment that $B$ is True or $A$ is False}

$$P(\bar A|T)=P(\bar A)/P(\bar A\mbox{\it\ or }B)
\mbox{\ \ \ and \ \ \ }
P(B|T)=P(B)/P(\bar A\mbox{\it\ or }B)$$

This case is the familiar identification of $T$ with $\bar A\mbox{\it\ or
}B$, the so-called material implication. The minimum amount of information
associated with the judgment is clearly $i(\bar A\mbox{\it\ or }B)$.

Unlike the other cases, the new probabilities of the propositions
$A$ and $B$, that is, the probabilities of those propositions after
the judgment, are not equal to old probabilities of logical expressions
such as $AB$ or $B\mbox{\it\ or }A$.

In this case, if $\bar T$ is given, $A\bar B$ is given, and
so the negation of the material implication, $\bar A\mbox{\it\ or }B$,
does not yield another relation of implication between the
relevant propositions.

The assignments of probabilities accomplished by this judgment
are left unchanged by the substitution of $\bar B$ for $A$ and
$\bar A$ for $B$, because the expression $\bar A\mbox{\it\ or }B$ is
left unchanged by this substitution. 

\vspace{5mm}

\subsection{The Judgment that $A$ is a Sufficient Condition of $B$}

\label{sufficientcondition}

If the proposition $B$ implies and is implied by the proposition
$B_1\mbox{\it\ or\ }B_2\mbox{\it\ or\ }\cdots B_N$ then
it is reasonable to call the propositions $B_1$, $B_2$ and so on the
sufficient conditions of $B$. 

If the judgment through which $T$ is given establishes $A$ as a new
sufficient condition of $B$, then $B$ will, after the judgment, 
imply and be implied by $B_1\mbox{\it\ or\ }\cdots B_N \mbox{\it\ or\ }A$.

This leads to:

$$P(A|T)=P(A)
\mbox{\ \ \ and \ \ \ }
P(B|T)=P(B\mbox{\it\ or }A)$$

In this case, the judgment that $A$ implies $B$ does not make $A$
any less likely, but it does make $B$ more likely. 

An example in English would be given by the sentence, ``If the weather
is good tomorrow, then I will go outside.'' For a protagonist who makes
this decision or a spectator who learns about it, the likelihood
that the weather will be good tomorrow does not increase, but the
likelihood that the protagonist will go outside does increase.

We can quantify the minimum amount of information which must be given
by the judgment by observing that $P(AB|T)=P(A|T)=P(A)$.  The amount of
information necessary to believe $AB$ therefore decreases from $i(AB)$
to $i(A)$ and so the amount of information provided, namely $i(T)$, must
be at least $i(AB)-i(A)$. This is equal to $i(B|A)$, which is perhaps
not surprisingly the quantity of information which it is necessary to
believe in order to believe that $B$ is true if $A$ can be assumed.

$i(T)$ is also related to the question of what happens if
$\bar T$ is given instead of $T$. Since the probability of $AB$
increases if $T$ is given, it can be expected to decrease if
$\bar T$ is given. It can be calculated explicitly as follows:
$$
\begin{array}{rl}
P(AB) & = P(ABT)+P(AB\bar T) \\
\Rightarrow P(AB\bar T) & =P(AB)-P(AB|T)P(T) \\
\Rightarrow P(AB\bar T) & =P(AB)-P(A)P(T) \\
\Rightarrow P(AB|\bar T) & =\frac{P(AB)-P(A)P(T)}{1-P(T)} \\
\end{array}
$$
which has two extreme solutions as $P(T)$ is varied, one of which is 
$P(AB|\bar T)=0$
at $P(T)=P(B|A)$ or $i(T)=i(B|A)$, and the other of which is
$P(AB|\bar T)=P(AB)$ at $P(T)=0$ or $i(T)=\infty$. The constraint
$i(T)\geq i(B|A)$ is necessary to ensure that all assigned
probabilities, such as $P(AB|\bar T)$, are non-negative numbers.
It applies only if $P(A)\neq 0$.

If $i(T)=i(B|A)$, then:
$$
i(AB;T)=\log\frac{P(ABT)}{P(AB)P(T)}=\log\frac{P(AB|T)}{P(AB)}
=\log\frac{P(A)}{P(AB)}=i(B|A)=i(T)
$$
which establishes that $AB$ implies $T$. This is the condition
under which the judgment that $T$ is true can accomplish
nothing independent of $AB$, and so in this sense, this is
exactly the amount of information required to establish
the implication and nothing else.

When $i(T)=i(B|A)$, $P(AB|\bar T)=0$, which shows that $A$
will imply $B$ if $T$ is given and $\bar B$ if $\bar T$ is
given.

Of the four cases discussed here, this form of judgment most closely
matches what is understood when the information provided by the simple
expression ``If $A$ then $B$'' is incorporated for the first time, since
that expression on its own does not indicate that the condition, $A$,
is any more or less likely to be satisfied. It does, however, indicate
a new condition under which $B$ can be known to be true, and $B$ then
becomes correspondingly more probable. 

After this section, when the expression ``If $A$ then $B$'' is used
to indicate a proposition, it should be understood that the
corresponding judgment establishes that $A$ is a sufficient
condition of $B$.

\vspace{5mm}

\subsection{The Judgment that $B$ is a Necessary Condition of $A$}

If the proposition $A$ implies and is implied by the proposition
$A_1A_2\cdots A_N$ then it is reasonable to call the propositions $A_1$,
$A_2$ and so on the necessary conditions of $A$.

If the judgment through which $T$ is given establishes $B$ as a new
necessary condition of $A$, then $A$ will, after the judgment, 
imply and be implied by $A_1\cdots A_N B$.

This leads to:

$$P(A|T)=P(AB)
\mbox{\ \ \ and \ \ \ }
P(B|T)=P(B)$$

Here, the judgment makes $A$ less likely without changing the
probability of $B$. In words it could be said that $B$ is judged
to be a necessary condition of $A$.

In English this is seen in the sentence, ``If the team is to win the
game, then they will need to score another goal.''  The probability
of the proposition that the team will win the game decreases when this
information is understood for the first time, while the probability of
the proposition that the team will score another goal is not increased
by the information provided in the sentence alone (although it might
increase if extra information is added, such as the information that
the team is trying to win the game).

If we replace $A$ with $\bar B$ and $B$ with $\bar A$ in the
probabilities above, we get:
$$
P(\bar B|T)=P(\bar A\bar B) \mbox{\ \ \ \ and\ \ \ \ } P(\bar A|T)=P(\bar A)
$$
which are equivalent to:
$$
P(B|T)=P(B\mbox{\it\ or }A) \mbox{\ \ \ \ and\ \ \ \ } P(A|T)=P(A)
$$
which is exactly the same as the previous case, in which $A$ was
judged to be a sufficient condition of $B$.

It can then be said that the judgment that $B$ is a necessary
condition of $A$ accomplishes the same as the judgment
that $\bar B$ is a sufficient condition of $\bar A$. This
allows us to say immediately that the information associated
with the judgment is $i(\bar A\bar B)-i(\bar B)$, or $i(\bar A|\bar B)$. 

However, it is evidently not the case that the judgment that $B$ is
a necessary condition of $A$ is equivalent to the judgment
that $A$ is a sufficient condition of $B$.

This illustrates the distinction between an assertion and a judgment. An
assertion merely involves a statement that a proposition is true, while
a judgment involves a decision to regard the proposition as true. A
judgment can change the circumstance, while an assertion is merely right
or wrong. Judgments and assertions therefore have different criteria of
equivalence: two judgments are equivalent if they effect the same change
in circumstance, while two assertions are equivalent if, whenever one
is right, the other is also right, and whenever one is wrong, the other
is also wrong.

The assertion that $B$ is a necessary condition of $A$ is then equivalent
to the assertion that $A$ is a sufficient condition of $B$, because
each proposition is true if and only if the other is. The respective
judgments of those propositions, however, are not equivalent, because one
judgment changes the conditions of, and hence the probability of, $A$,
while the other changes the conditions of $B$.

\vspace{5mm}

\subsection{\small{Judging that $A$ Implies $B$ Without Changing Either $P(A)$ or $P(B)$}}

$$P(A|T)=P(A) 
\mbox{\ \ \ and \ \ \ }
P(B|T)=P(B)$$

In addition, since the judgment establishes that $A$ implies $B$,
it is also true that $P(AB|T)=P(A|T)=P(A)$ and $P(\bar A\bar B|T)=P(\bar
B|T)=P(\bar B)$.

If the judgment is to leave the probabilities of $A$ and $B$ unchanged
but is to establish that $A$ implies $B$ then it must be the case before
the judgment that the probability of $B$ is greater than or equal to
the probability of $A$, or, equivalently, that the amount of information
which it is necessary to believe in order to believe $A$ is greater
than or equal to the amount of information necessary to believe $B$.
This kind of judgment can therefore only be made in certain
circumstances, unlike the other three. 

In a sense, the judgment places the information which it is necessary
to believe in order to believe $B$ within the information which it
is necessary to believe in order to believe $A$. At the same time,
it places the information required to believe $\bar A$ within the
information required to believe $\bar B$. Like material implication,
the effect of the judgment is invariant if $A$ is replaced by $\bar B$
and $B$ is replaced by $\bar A$.

To quantify the amount of information associated with the judgment,
we can observe that:
$$
\begin{array}{rl}
P(AB) & = P(ABT)+P(AB\bar T) \\
\Rightarrow P(AB|\bar T) & =\frac{P(AB)-P(A)P(T)}{1-P(T)} \\
\end{array}
$$
as in the second case, which implies that $i(T)\geq i(B|A)$.

In this case, however, there is the additional constraint that
$$
\begin{array}{rl}
P(\bar A\bar B) & = P(\bar A\bar BT)+P(\bar A\bar B\bar T) \\
\Rightarrow P(\bar A\bar B|\bar T) & =\frac{P(\bar A\bar B)-P(\bar B)P(T)}{1-P(T)} \\
\end{array}
$$
which gives the constraint $i(T)\geq i(\bar A|\bar B)$. 

In order for both of these to be satisfied, $i(T)$ must be
greater than or equal to whichever of $i(B|A)$ and
$i(\bar A|\bar B)$ is greater. 

We can also see that $i(AB;T)=i(B|A)$ and $i(\bar A\bar B;T)=i(\bar
A|\bar B)$ using the same reasoning which led to $i(AB;T)=i(B|A)$
in the second case. Taking $i(T)$ to equal whichever of $i(B|A)$ and
$i(\bar A|\bar B)$ is greater, then, we see that either 
$i(AB;T)=i(T)$ or $i(\bar A\bar B;T)=i(T)$, so this judgment
would establish the implication and nothing else, meaning
nothing independent of $AB$ if $i(B|A)\geq i(\bar A|\bar B)$ or
nothing independent of $\bar A\bar B$ if $i(\bar A|\bar B)\geq i(B|A)$.

In the case when $i(B|A)\geq i(\bar A|\bar B)$ and $i(T)=i(B|A)$,
it can easily be seen from the above that $P(AB|\bar T)=0$, so
if $\bar T$ is given then $A$ will imply $\bar B$ and $B$ will
imply $\bar A$. 

If, instead, $i(\bar A|\bar B)\geq i(B|A)$ and $i(T)=i(\bar A|\bar B)$,
then $P(\bar A\bar B|\bar T)=0$, so if $\bar T$ is given then
$\bar A$ will imply $B$ and $\bar B$ will imply $A$.

\subsection{Counterfactuals}

If the material conditional, $B\mbox{\it\ or }\bar A$, is the only way
which can be used to represent ``$A$ implies $B$'', then it must
be said that a true proposition is implied by any proposition
and a false proposition implies every proposition. 

If, however, ``$A$ implies $B$'' is represented by $P(B|A)=1$, then
it is not possible for $A$ to imply both $B$ and $\bar B$. If
$P(B|A)=1$, then $P\left(\bar B|A\right)$ must equal zero because
$P(B|A)$ is the probability that $B$ would have if $A$ were 
to be given. If $A$ is given, $P(B)+P(\bar B)$ will need
to equal one, regardless of what the probability of $A$ was
before the judgment.

Reflecting on the fact that a proposition, $A$, with $P(A)=0$ in the
current circumstance, does not already imply every other proposition,
$B$, in the sense that $P(B|A)=1$, we can ask how much information is
given when it is judged that $A$ implies $B$. That is, we can
consider what is involved in a judgment which changes the
value of $P(B|A)$ to one.

In particular, we can consider the proposition, $T$, which establishes
that $A$ is a sufficient condition of $B$, so that $P(B|T)=P(B\mbox{\it\
or\ }A)=P(B)$ and $P(A|T)=P(A)=0$.  This evidently coincides, in this
case, with the judgment which establishes the implication without changing
either $P(A)$ or $P(B)$.

In section \ref{sufficientcondition}, it was shown that,
when $P(A)\neq 0$, $i(T)$ must be greater than or equal
to $i(B|A)$. The proof relied on the fact that $i(AB)$
decreases to $i(A)$ when $T$ is judged to be true, and
so $i(T)$ must account for this decrease, $i(AB)-i(A)$.
If $P(A)=0$, though, the expression $i(AB)-i(A)$ gives
$\infty-\infty$. 

The other way of proving that $i(T)\geq i(B|A)$ relied on the fact that:
$$ 
P(AB|\bar T) =\frac{P(AB)-P(A)P(T)}{1-P(T)} 
$$
and so $P(T)$ must be less than or equal to $P(AB)/P(A)$ for
$P(AB|\bar T)$ to be non-negative. If $P(A)=0$, however, this
does not apply, because $P(AB)-P(A)P(T)=0$, and so the lower bound on
$i(T)$, namely $i(B|A)$, does not apply if $P(A)=0$.

There is, however, a different condition which must be
satisfied, because:
$$
P(B|A)=P(BT|A)+P(B\bar T|A)=P(B|TA)P(T|A)+P(B|\bar TA)P(\bar T|A)
$$
In combination with $P(B|TA)=1$, this establishes that:
$$
P(B|A)\geq P(T|A)
$$
or
$$
i(T|A)\geq i(B|A)
$$

When $T$ is the proposition that $A$ is a sufficient condition of $B$,
in a circumstance where $P(A)=0$, then, the condition $i(T|A)\geq i(B|A)$
applies but $i(T)\geq i(B|A)$ does not. 

\subsubsection{When a Proposition with Probability One is Given}

In fact, there is no reason why $i(T)$ can not equal zero, or,
equivalently, why $P(T)$ can not equal one. If $P(A)$ were
greater than zero, $i(T)=0$ would imply $i(T|A)=0$.
However, if $P(A)=0$, then $P(T|A)$ is not constrained
by the values of $P(A)$ and $P(TA)$, so $i(T|A)$ can
take any value, even if $i(T)=0$. 

When $P(T)=1$, it is guaranteed that $P(X|T)=P(X)$ for every proposition,
$X$, because:
$$
P(X)=P(XT)+P(X\bar T)=P(X|T)P(T)+P(X|\bar T)P(\bar T)=P(X|T)
$$
It is not, however, guaranteed that $P(X|TA)=P(X|A)$. In
order to guarantee $P(X|TA)=P(X|A)$, the relevant condition
is $P(T|A)=1$, not $P(T)=1$.

There is, then, the peculiar possibility that a judgment which gives
a proposition, $T$, with $P(T)=1$, can actually produce a change in
the value of $P(B|A)$, although it cannot change the value of $P(X)$
for any proposition, $X$.

This illustrates the distinction between {\it believing} a proposition
and {\it being given} the proposition. As a more explicit example,
if the proposition ``$B\mbox{\it\ or not-}A$'' is believed, then
$$
\begin{array}{rl}
P(B\mbox{\it\ or }\bar A)=1 & \Leftrightarrow P(B)+P(\bar A)-P(B\bar A)=1 \\
& \Leftrightarrow P(B)+1-P(A)- \left(P(B)-P(AB)\right)=1 \\
& \Leftrightarrow 1-P(A)-P(B)+P(AB)=1-P(B) \\
& \Leftrightarrow P(AB)=P(A) 
\mbox{\ or equivalently\ } P\left(\bar A\bar B\right)=P\left(\bar B\right)\\
& \Leftrightarrow P(B|A)P(A)=P(A) 
\mbox{\ or equivalently\ } P\left(\bar A|\bar B\right)P(\bar B)=P\left(\bar B\right)\\
\end{array}
$$
so either:
\begin{itemize}
\item
$B$ is believed, $P\left(\bar B\right)=0$, or
\item
$\mbox{not-}A$ is believed, $P(A)=0$, or
\item
$P(B|A)=1$ and $P\left(\bar A|\bar B\right)=1$.
\end{itemize}

That is, if ``$B\mbox{\it\ or not-}A$'' is believed, then either
$B$ is believed or $\mbox{not-}A$ is believed or the truth of
each can be inferred from the falsity of the other.

``$B\mbox{\it\ or not-}A$'' can therefore be believed without it
being possible to infer $B$ from $A$. If $P(A)=0$, then $P(B|A)$
can take any value, even 0, while $P\left(B\mbox{\it\ or }\bar A\right)=1$.

On the other hand, if the proposition ``$B\mbox{\it\ or not-}A$'' is
given, rather than merely believed, then it can be used in
combination with $A$ to infer $B$, so $P(B|A)=1$, and, similarly,
$P\left(\bar A|\bar B\right)=1$.

A proposition with probability one, therefore, such as ``$B\mbox{\it\
or not-}A$'' when $P(A)=0$, can change the values of hypothetical, or
conditional, probabilities, such as $P(B|A)$, when it is given, although
it can not change the values of any actual probabilities.  

It can be observed that, when $P(A)=0$, the changes produced when
``$B\mbox{\it\ or not-}A$'' is given are identical to those produced
by the judgment that $A$ is a sufficent condition of $B$. In
each case, $P(B|A)$ changes to one, while the probabilities of
propositions are unchanged.

It may seem inappropriate to say that no information is provided when
a proposition with probability one is given, since the judgment
actually produces a change. 

The changes which can be made by a judgment through which a proposition
with probability one is given are always changes to quantities such as
$P(B|A)$ where $A$ is a proposition of probability zero. They will never
have any influence on the probabilities of propositions, such as $B$,
unless a circumstance is reached in which $P(A)\neq 0$. In order for
this to happen, a proposition of probability zero (for example, $A$)
must be given.

It can then be said that no information is provided when a proposition,
$T$, of probability one is given, until a proposition of probability
zero, such as $A$, is given later, at which point the amount $i(T|A)$
of counterfactual information associated with the earlier judgment
becomes relevant.  This exception always applies, though.  The judgment
that $A$ is true, for example, establishes only that $A$ will be
believed until a proposition of probability zero is given. The caveat
``until a proposition of probability zero is given'' applies to any
description of what is accomplished by a judgment.

\section{The Paradox of the Raven}
\label{raven}

\subsection{The Judgment that all Ravens are Black}

The paradox \cite{hempel1} can be summarized as follows: 
\begin{itemize}
\item
The statement that all ravens
are black is logically equivalent to the statement that all non-black
things are non-ravens. 
\item
An observation which provides evidence supporting
one statement, such as the observation of a black raven, therefore
also supports the other statement. 
\item
Therefore the observation of a
non-black non-raven, such as a blue sky, provides evidence in
favour of the statement that all ravens are black.
\end{itemize}

The universal proposition, ``All ravens are black,'' can be expressed in
the form ``If $x$ is a raven then $x$ is black,'' where $x$ is understood
to be a variable. Any name can be substituted for $x$ in the expression
above to form a singular proposition which is implied by the universal
one, for example, ``If Socrates is a raven then Socrates is black.''

Now ``Socrates is a raven'' and ``Socrates is black'' are propositions
which take the place of the propositions $A$ and $B$ from the previous
section, where it was shown that ``If $A$ then $B$'' can be interpreted
in at least four ways. One can ask which of the four ways is
understood when it is judged that all ravens are black on the basis
of evidence in the form of observations of black ravens. 

The first question is whether, in making such a judgment, one
intends to reduce the probability that ``Socrates is a raven,''
since two of the four ways of judging that ``If $A$ then $B$''
reduce the probability of $A$, namely the judgment that either
$B$ is true or $A$ is false and the judgment that $B$ should
from now on be regarded as a necessary condition of $A$.

Upon inspection, the judgment that all ravens are black does not
appear to involve any reduction in the probability that something
is a raven. After the judgment, one does not think that there
are fewer ravens than one had previously thought there were. One
merely thinks that the ravens there are, however many there may be,
are black, where it was previously uncertain what colour they were. 

The two remaining ways of judging that ``If $A$ then $B$'' are the
judgment that $A$ is a sufficient condition of $B$, which sets the new
value of $P(B)$ to the old value of $P(B\mbox{\it\ or\ }A)$ without
changing the value of $P(A)$, and the judgment which establishes the
implication without changing either $P(A)$ or $P(B)$. 

The second question, then, is whether the judgment that all ravens
are black, made after the observation of many black ravens, changes
the probability that ``Socrates is black'' by setting it equal to
the old value of ``Socrates is black {\it or} Socrates is a raven'',
or leaves it unchanged. 

If it was the case before the judgment that Socrates was believed to be a
raven but it was uncertain whether he was black, then after the judgment
he is believed to be black. The overall effect is then to increase the
probability that ``Socrates is black'' by setting it equal to the old
value of ``Socrates is black {\it or} Socrates is a raven.''

The judgment that all ravens are black, made after observing black ravens,
therefore makes ``$x$ is a raven'' a sufficient condition of ``$x$
is black'', thereby increasing the probability that $x$ is black without
changing the probability that $x$ is a raven. 

Correspondingly, the judgment that all non-black things are non-ravens
would increase the probability that ``$x$ is a non-raven'' and leave
the probability that ``$x$ is not black'' unchanged. It would therefore
leave the probability that $x$ is black unchanged while decreasing the
probability that $x$ is a raven. 

\subsection{\small{The Evidence Provided by the Observation of a Black Raven}}

The different judgments, ``All ravens are black,'' and ``All non-black
things are non-ravens,'' require different evidence. In both cases,
it is possible to quantify the amount of evidence provided by the 
observation of a black raven.

Let $A$ be the proposition that ``All ravens are
black''\footnote{Whether this refers to all logically possible ravens
or merely all the ravens in a particular collection is not addressed
here. In the former case, $P(A)$ would have to equal zero (since we are
supposing that ravens are not by definition black) and in the latter
case $P(A)$ could be greater than zero.}, and let $N$ be the
proposition that ``All non-black things are non-ravens.'' Let $R$ be
the proposition that ``Socrates is a raven'' and let $B$ be the
proposition that ``Socrates is black.''

The amount of evidence provided by the proposition that Socrates
is black, in favour of the proposition that all ravens are
black, given that Socrates is a raven, is\footnote{We are considering
a circumstance in which $R$ is actually given, but we occasionally
include $R$ in expressions such as $P(B|AR)$ to serve as
a reminder of this fact.}:
$$
e(B\rightarrow A|R)=\log\frac{P(B|AR)}{P\left(B|\bar AR\right)}
$$

It is clear that ``All ravens are black'' implies that
``Socrates is black'' when ``Socrates is a raven'' is given,
so:
$$
e(B\rightarrow A|R)=\log\frac{1}{P\left(B|\bar AR\right)}=i\left(B|\bar AR\right)
$$

Since it is not already believed that Socrates is black, $1>P(B|R)$,
which, together with:
$$
P(B|R)=P(B|AR)P(A|R)+P\left(B|\bar AR\right)P\left(\bar A|R\right)
$$ 
and $P(B|AR)=1$, implies that $P\left(B|\bar AR\right)<1$, leading to:
$$
e(B\rightarrow A|R)=i\left(B|\bar AR\right)>0
$$
which proves that discovering that ``Socrates is black'' when
it is already believed that ``Socrates is a raven'' does in fact
provide evidence in favour of the proposition that ``All ravens are black''.

If we want to calculate the amount of evidence provided
by the proposition that Socrates is black in favour of the proposition
that all non-black things are non-ravens, we must first calculate:
$$
e(B\rightarrow N|R)=\log\frac{P(B|NR)}{P\left(B|\bar NR\right)}
$$
which involves calculating the value of $P(B|NR)$ which is the
probability that $B$ would have if $N$ were to be judged to
be true after $R$ had already been given. 

Now $B$ is the proposition that Socrates is black and $N$ is the
proposition that all non-black things are non-ravens.  $N$ clearly implies
``If Socrates is not black then Socrates is not a raven,'', or ``If $\bar
B$ then $\bar R$.''  If $N$ is judged to be true, the new probability that
``Socrates is a non-raven'' will be equal to the old probability that
``Socrates is a non-raven {\it or } Socrates is not black.''

If $N$ is judged to be true, then it is not necessarily believed that
``Socrates is a raven'' any more, since the later judgment that ``If
Socrates is not black then Socrates is not a raven,'' conditionally
overrules the earlier judgment that ``Socrates is a raven''. This
illustrates another distinction between judgments and
assertions, namely that assertions can contradict other assertions,
while judgments overrule earlier judgments.

Overruling the earlier judgment requires an infinite amount of
information: it was previously completely certain that ``Socrates
is a raven'', $P(R)=1$ or $P\left(\bar R\right)=0$ or $i\left(\bar
R\right)=\infty$, but $i\left(\bar R|N\right)=-\log P\left(\bar
R\mbox{\it\ or }\bar B\right)<\infty$. When $N$ is given, the
amount of information which it necessary to believe in order to
believe that $R$ is false decreases from infinity to a finite value,
and a correspondingly infinite amount of information must be given. That
is, $P(N)=0$.

This can be confirmed by observing that the amount of information
required to establish ``If $\bar B$ then $\bar R$'' is $i\left(\bar R|\bar
B\right)$, which is infinite if $R$ is already believed to be true. Since
$N$ implies ``If $\bar B$ then $\bar R$'', $i(N)\geq i\left(\bar R|\bar
B\right)=\infty$.

The proposition that ``All non-black things are non-ravens'' predicates
non-ravenness of various things, including Socrates if he should turn
out to be a non-black thing, but it does not predicate blackness of
anything. The probability that Socrates is black, $P(B)$, is therefore
unaffected by the judgment, just as the judgment that ``If $\bar B$
then $\bar R$'' does not change the probability of $\bar B$.

The result is that $P(B|NR)=P(B|R)$, and it follows from $P(N)=0$
that $P\left(B|\bar NR\right)=P(B|R)$, so:
$$
e(B\rightarrow N|R)=\log\frac{P(B|R)}{P\left(B|R\right)}=0
$$

The result does not change if it is initially believed that ``Socrates
is black'' and then subsequently discovered that ``Socrates is a
raven'': In circumstances where ``Socrates is black'' is believed, the
proposition ``If Socrates is not black then Socrates is not a raven'' is
a counterfactual so when it is given it does not change the probability
that ``Socrates is a raven''.

The judgment that ``All non-black things are non-ravens'' therefore does
not change the probability that ``Socrates is a raven'', $P(R|NB)=P(R|B)$.  
This implies that $P\left(R|\bar NB\right)=P(R|B)$ since 
$P(R|B)=P(R|NB)P(N|B)+P\left(R|\bar NB\right)P\left(\bar N|B\right)$
and $P\left(\bar N|B\right)\neq 0$, leading to: 
$$
e(R\rightarrow N|B)=\log\frac{P(R|B)}{P\left(R|B\right)}=0
$$

This shows that the amount of evidence, provided by the
observation of a black raven, in favour of the proposition
that all non-black things are non-ravens, is zero. Correspondingly,
the observation of a non-black non-raven, such as a blue sky,
provides zero evidence in favour of the proposition that all
ravens are black.

This resolves the paradox.

\section*{{\small Acknowledgements}}
{\small The author would like to thank Charles F. Stevens, Joe Snider,
Posina Venkata Rayudu, Doug Rubino, Bayle Shanks and Jennifer
Bowen for useful conversations and comments on the manuscript.}

\vspace{5mm}

\bibliographystyle{prsty}
\bibliography{judgmentbib}

\end{document}